\documentclass{amsart}

\usepackage[utf8]{inputenc}

\usepackage{color}
\usepackage[colorlinks=true]{hyperref}
\hypersetup{urlcolor=blue, citecolor=blue}

\usepackage{amsmath,amssymb,amsfonts,amsthm}
\usepackage{graphicx}
\usepackage{pgfplots}
\usepackage{cleveref}
\usepackage{thmtools}
\usepackage{stmaryrd}
\usepackage{mathtools}
\usepackage{mathrsfs}  
\usepackage{enumerate}
\usepackage{dsfont}

\usepackage[normalem]{ulem} 
\usepackage{cancel}
\usepackage[color=red!65!yellow]{todonotes}

\usepackage{subfig}
\usepackage{booktabs}

\usepackage[ruled,vlined]{algorithm2e}

\begin{document}

\title{TPFA Finite Volume approximation of Wasserstein gradient flows}

\author[A. Natale]{Andrea Natale}
\address{Andrea Natale (\href{mailto:andrea.natale@u-psud.fr}{\tt andrea.natale@u-psud.fr}) Laboratoire de Mathematiques d'Orsay, Universit\'e Paris-Sud
} 
\author[G. Todeschi]{Gabriele Todeschi}
\address{Gabriele Todeschi (\href{mailto:gabriele.todeschi@inria.fr}{\tt gabriele.todeschi@inria.fr}):  Inria, Project team Mokaplan, Universit\'e Paris-Dauphine, PSL Research University, UMR CNRS 7534-Ceremade
}

\maketitle

\begin{abstract}
Numerous infinite dimensional dynamical systems arising in different fields have been shown to exhibit a gradient flow structure in the Wasserstein space. We construct Two Point Flux Approximation Finite Volume schemes discretizing such problems which preserve the variational structure and have second order accuracy in space. We propose an interior point method to solve the discrete variational problem,  providing an efficient and robust algorithm.
We present two applications to test the scheme and show its order of convergence.
\end{abstract}

\section{Gradient flows' time discretization}

A gradient flow is a process that, starting from an initial point, evolves by maximizing at each instant the rate of decay of a given specific energy. Many problems arising in physics, biology, social sciences, etc., can be recast as infinite dimensional gradient flows. Considering a compact domain $\Omega\subset\mathbb{R}^d$, a finite time horizon $T\in\mathbb{R}^+$, and a real-valued, strictly convex and proper energy functional $\mathcal{E}$, we focus our attention on problems of the form
	\begin{equation}\label{eq:pde}
	\begin{cases}
	\partial_t \rho - \nabla \cdot (\rho \nabla \frac{\delta \mathcal{E}}{\delta \rho}[\rho]) = 0, \quad &\text{in} \; \Omega\times[0,T], \\
	\rho \nabla \frac{\delta \mathcal{E}}{\delta \rho} \cdot \boldsymbol{n} = 0, \quad &\text{on} \; \partial\Omega\times[0,T], \\
	\rho(0) = \rho^0, \quad &\text{in} \; \Omega,
	\end{cases}
	\end{equation}
	where $\frac{\delta \mathcal{E}}{\delta \rho}$ denotes the first variation of $\mathcal{E}$, $\rho^0\in L^1(\Omega;\mathbb{R}^+)$ is a given initial condition and $\boldsymbol{n}$ is the unit outer normal vector to $\partial \Omega$. Problem \eqref{eq:pde} denotes the continuity equation of a time evolving non-negative density $\rho$ convected by the velocity field $- \nabla \frac{\delta \mathcal{E}}{\delta \rho}[\rho]$, with no flux across the boundary of the domain, hence preserving its total mass.
	It is nowadays clear that problems of the form of \eqref{eq:pde} represent gradient flows of the energy $\mathcal{E}$ with respect to the Wasserstein metric.
	We refer to \cite{AGS08,Santambrogio_OTAM} for more details on gradient flows and optimal transport.
	
	The underlying variational structure of this type of problems provides useful tools for their study. From the numerical point of view, more robust solvers can be designed by taking it into account. In particular, the property that the energy should decrease as fast as possible at each time step is a useful criterion to assess the goodness and reliability of a numerical solution and it should be preserved.
	The JKO scheme realizes this by using the variational formulation of the implicit Euler method.
	For an increasing sequence $(t^n)_{n\in\mathbb{N}}\subset\mathbb{R}$ of time steps such that $\cup_n [t^{n-1},t^n]=[0,T]$, let $Q^n=\Omega\times[t^{n-1},t^n]$ and $\partial Q^n=\partial\Omega\times[t^{n-1},t^n]$. The JKO scheme constructs a sequence $(\rho^n)_{n\in\mathbb{N}}$ as follows: given an approximation $\rho^{n-1}$ of the density at time $t^{n-1}$, compute $\rho^n=\tilde{\rho}(t^n)$, where $(\tilde{\rho},\tilde{\boldsymbol{F}}):Q^n \rightarrow \mathbb{R}^+\times\mathbb{R}^d$ solve
	\begin{equation}\label{eq:jko}
	\inf_{(\tilde{\rho},\tilde{\boldsymbol{F}})} \int_{Q^n} \frac{|\tilde{\boldsymbol{F}}|^2}{2\tilde{\rho}} \text{d} \boldsymbol{x} \text{dt}+ \mathcal{E}(\tilde{\rho}(t^n)),\; \;\text{where $(\tilde{\rho},\tilde{\boldsymbol{F}})$ solve:}\; 
	\begin{cases}
	\partial_t \tilde{\rho} + \nabla \cdot \tilde{\boldsymbol{F}} = 0,  &\text{in} \; Q^n\\
	\tilde{\boldsymbol{F}}\cdot \boldsymbol{n} = 0, &\text{on} \;  \partial Q^n\\
	\tilde{\rho}(t^{n-1}) = \rho^{n-1}.\\
	\end{cases}
	\end{equation}
	The density $\rho^n$ is computed minimizing the sum of its squared Wasserstein distance from $\rho^{n-1}$ and the energy in $\rho^n$. The former term corresponds to the total kinetic energy of the curve $\tilde{\rho}$ written in the variables density-momentum, $(\tilde{\rho},\tilde{\boldsymbol{F}})$, rather than density-velocity, in order to highlight the convexity of the problem \cite{BB00}.
	The sequence of densities $(\rho^n)_{n\in\mathbb{N}}$, meant to be an approximation of the solution at each time step $t^n$, can be seen as a piecewise constant time-dependent density converging to the flow under suitable assumptions \cite{AGS08,Santambrogio_OTAM}. This time discretization enables to design energy-diminishing schemes that are furthermore robust in the sense that, since \eqref{eq:jko} is a well-posed convex problem, the solution at step $n$ always exists no matter the time step $\tau^n =t^n-t^{n-1}$.

The Wasserstein distance involved in \eqref{eq:jko} needs to be further discretized in time. Since the JKO scheme is of order one \cite{JKO98}, a first order time discretization is sufficient and leads to a reasonable computational complexity. 
We can approximate \eqref{eq:jko} with an LJKO \cite{LJKO_arXiv}: given an approximation $\rho^{n-1}$ of the density at time $t^{n-1}$, compute $\rho^n$ solution to
	\begin{equation}\label{eq:ljko}
	\inf_{(\rho,\boldsymbol{F})} \tau^n \int_{\Omega} \frac{|\boldsymbol{F}|^2}{2\rho} \text{d} \boldsymbol{x} + \mathcal{E}(\rho), \;\; \text{where $(\rho,\boldsymbol{F})$ solve:} \;
	\begin{cases}
	\rho-\rho^{n-1} + \tau^n \; \nabla \cdot \boldsymbol{F} = 0,  &\text{in} \; \Omega,\\
	\boldsymbol{F} \cdot \boldsymbol{n} = 0, &\text{on} \; \partial\Omega,
	\end{cases}
	\end{equation}
	where now $(\rho,\boldsymbol{F}):\Omega\rightarrow\mathbb{R}^+\times\mathbb{R}^d$ does not depend on time. The continuity equation is discretized using a single implicit Euler step, whereas the time integral using a right endpoint approximation.

Given the conservative form of the problem, Finite Volume methods appear as natural choices for its discretization. Their relation with optimal transport has been highlighted in, e.g., \cite{GKM_arXiv}.
	Ensuring the positivity of the density is a crucial property for any candidate numerical method, since problems \eqref{eq:jko} and \eqref{eq:ljko} lose their convexity if the density is negative.
	In \cite{LJKO_arXiv} problem \eqref{eq:ljko} is discretized using upwind FV, which provides automatically the positivity for the discrete solution. The problem can then be solved using a Newton scheme. However, this gives an order one space discretization. Moreover, the derived scheme is not particularly robust since small time steps may be required to make the Newton scheme converge.
	In the present work we propose a more general FV framework, which allows us to consider second order discretizations in space. As a consequence, the positivity constraint on the density needs to be taken into account. To this end, we use an interior point method.	
	
\section{Finite Volume discretization}
	
	Assume the domain $\Omega \subset \mathbb{R}^d$ to be polygonal if $d=2$ or polyhedral if $d=3$. The specifications for a partitioning of $\Omega$ to be admissible for TPFA Finite Volume are classical \cite[Definition 9.1]{EGH00}. We denote by $\left(\mathcal{T}, \overline{\Sigma}, {(\boldsymbol{x}_K)}_{K\in\mathcal{T}}\right)$ such an admissible mesh, namely the triplet of the set of polyhedral control volumes, the set of faces and the set of cell centers. We use Delaunay triangulations in order to satisfy these assumptions. The Lebesgue measure of $K\in\mathcal{T}$ is denoted by $m_K>0$.
	The set $\overline{\Sigma}$ is composed of boundary faces $\Sigma_{ext} = \{ \sigma \subset \partial\Omega\}$ and internal faces $\sigma \in \Sigma = \overline{\Sigma} \setminus \Sigma_{ext}$. We denote by $\Sigma_{K} = \overline{\Sigma}_{K}\cap \Sigma$ the internal faces belonging to $\partial K$. For each internal face $\sigma = K|L \in \Sigma$, we refer to the diamond cell $\Delta_\sigma$ as the polyhedron whose edges 
	join $\boldsymbol{x}_K$ and $\boldsymbol{x}_L$ to the vertices of $\sigma$. Denoting by $m_{\sigma}$ the Lebesgue measure of the edge $\sigma$ and by $d_\sigma = |\boldsymbol{x}_K-\boldsymbol{x}_L|$, the measure $m_{\Delta_{\sigma}}$ of $\Delta_\sigma$ is then equal to $m_\sigma d_\sigma/d$, where 
	$d$ stands for the space dimension. We denote by $d_{K,\sigma}$ the euclidean distance between the cell center $\boldsymbol{x}_K$ and the midpoint of the edge $\sigma \in \overline{\Sigma}_K$. The size of the mesh is defined by $h_{\mathcal{T}}=\max_{K\in\mathcal{T}} \text{diam}(K)$.
	
	We introduce the space of discrete conservative fluxes
	\[
	\mathbb{F}_{\mathcal{T}}=\{\boldsymbol{F}=(F_{K,\sigma},F_{L,\sigma})_{\sigma\in\Sigma}\in\mathbb{R}^{2\Sigma}: F_{K,\sigma}+F_{L,\sigma}=0\}
	\]
	and denote $F_{\sigma} = |F_{K,\sigma}| = |F_{L,\sigma}|$. We introduce also the spaces of discrete variables on cells $\mathbb{P}_{\mathcal{T}} = \mathbb{R}^{\mathcal{T}}$ and diamond cells $\mathbb{P}_{\Sigma} = \mathbb{R}^{\Sigma}$, endowed with the two scalar products $\langle \cdot, \cdot \rangle_K: (\boldsymbol{a},\boldsymbol{b})\in [\mathbb{P}_{\mathcal{T}}]^2 \mapsto \sum_{K\in\mathcal{T}} a_K b_K m_K$, $\langle \cdot, \cdot \rangle_{\sigma}: (\boldsymbol{u},\boldsymbol{v})\in [\mathbb{P}_{\Sigma}]^2 \mapsto \sum_{\sigma\in\Sigma} u_{\sigma} v_{\sigma} m_{\sigma} d_{\sigma}$, respectively. We introduce a reconstruction operator from cells to diamond cells $R_{\Sigma}:\mathbb{P}_{\mathcal{T}}\rightarrow \mathbb{P}_{\Sigma}$. On each edge $\sigma = K|L$, the density on the diamond cell can be reconstructed from the values of the densities $\rho_K,\rho_L$.
	To keep the scheme simple, we employ weighted arithmetic averages $\rho_{\sigma}= \lambda_{K,\sigma} \rho_K + \lambda_{L,\sigma} \rho_L$, with $\lambda_{K,\sigma}, \lambda_{L,\sigma} \in [0,1], \lambda_{K,\sigma}+\lambda_{L,\sigma} = 1$. Nonetheless, other choices are possible, such as geometric, harmonic and logarithmic averages and all their weighted versions \cite{GKM_arXiv}.
	We consider three possibilities for the weights  $(\lambda_{K,\sigma}, \lambda_{L,\sigma})$: $(\frac{1}{2},\frac{1}{2})$, the standard arithmetic mean; $(\frac{d_{L,\sigma}}{d_{\sigma}},\frac{d_{K,\sigma}}{d_{\sigma}})$, which provides a linear reconstruction of the density at the edge midpoint; $(\frac{d_{K,\sigma}}{d_{\sigma}},\frac{d_{L,\sigma}}{d_{\sigma}})$, which gives a mass weighted arithmetic mean.
	Thanks to these choices we expect to obtain second order accuracy for the space discretization.
	We introduce also the adjoint operator of this reconstruction, with respect to the two scalar products, given by $R_{\mathcal{T}}:\boldsymbol{\rho}\in\mathbb{P}_{\Sigma} \mapsto \big( \sum_{\sigma\in\Sigma_K} \rho_{\sigma} \lambda_{K,\sigma} \frac{m_\sigma d_{\sigma}}{m_K}\big)_{K\in\mathcal{T}}\in \mathbb{P}_{\mathcal{T}}$.

	Assuming the energy $\mathcal{E}(\rho)$ to be of the form $\int_{\Omega} E(\rho) \text{d} \boldsymbol{x}$ for a real valued and strictly convex scalar function $E$, given the discrete initial density of the form $(\rho^0_K)_{K\in\mathcal{T}}=(\rho^0(\boldsymbol{x}_K))_{K\in\mathcal{T}} \in \mathbb{P}^+_{\mathcal{T}}$, the discrete LJKO scheme is: given $\boldsymbol{\rho}^{n-1}=(\rho^{n-1}_K)_{K\in\mathcal{T}}\in\mathbb{P}^+_{\mathcal{T}}$ approximation of the density at time $t^{n-1}$, compute $\boldsymbol{\rho}^n$ solution to
	\begin{equation}\label{eq:ljko_h}
	\inf_{(\boldsymbol{\rho},\boldsymbol{F})} \tau^n \sum_{\sigma\in\Sigma} \frac{F_{\sigma}^2}{2(R_{\Sigma}(\boldsymbol{\rho}))_{\sigma}} m_{\sigma}d_{\sigma} + \sum_{K\in\mathcal{T}} E(\rho_K) m_K,
	\end{equation}
	with $(\boldsymbol{\rho},\boldsymbol{F})\in \mathbb{P}_{\mathcal{T}}\times\mathbb{F}_{\mathcal{T}}$ such that $(\rho_K-\rho^{n-1}_K) m_K + \tau^n \;\sum_{\sigma\in\Sigma_K} F_{K,\sigma} m_{\sigma} = 0$ and $\rho_K\ge0, \forall K \in \mathcal{T}$. We take as measure of the diamond cell $d m_{\Delta_{\sigma}}$, as it is classically done in order to compensate the unidirectional discretization of the momentum \cite{EGH00}.
	The constraint $\boldsymbol{F}\cdot\boldsymbol{n}=0$ is automatically taken into account disregarding the flux on the boundary edges in the definition of the space of discrete conservative fluxes. The conservation of mass is also automatically enforced thanks to the conservativity of the Finite Volume discretization, i.e. $\sum_{K\in\mathcal{T}} \rho^n_K m_K = \sum_{K\in\mathcal{T}} \rho^{n-1}_K m_K$. Furthermore, the scheme guarantees a discrete energy-dissipation property: given the couple $(\boldsymbol{\rho}^n,\boldsymbol{F}^n)$ solution to \eqref{eq:ljko_h}, the competitor $(\boldsymbol{\rho}^{n-1},\boldsymbol{0})$ provides
	\[
	\tau^n \sum_{\sigma\in\Sigma} \frac{(F^n_{\sigma})^2}{2(R_{\Sigma}(\boldsymbol{\rho}^n))_{\sigma}} m_{\sigma}d_{\sigma} + \sum_{K\in\mathcal{T}} E(\rho^n_K) m_K \le \sum_{K\in\mathcal{T}} E(\rho^{n-1}_K) m_K.
	\]
	
	At each step $n$, \eqref{eq:ljko_h} is a strictly convex optimization problem with linear constraints. Enforcing the constraints with the multipliers $-\boldsymbol{\phi}\in\mathbb{P}_{\mathcal{T}},\boldsymbol{\lambda}\in\mathbb{P}^-_{\mathcal{T}}$ and using the definition of the conservative fluxes we obtain the saddle point problem
	\begin{multline}\label{eq:infsup_h}
	\inf_{(\boldsymbol{\rho},\boldsymbol{F})} \sup_{(\boldsymbol{\phi},\boldsymbol{\lambda})} \tau^n \sum_{\sigma\in\Sigma} \frac{(F_{\sigma})^2}{2(R_{\Sigma}(\boldsymbol{\rho}))_{\sigma}} m_{\sigma}d_{\sigma} + \sum_{K\in\mathcal{T}} (\rho^{n-1}_K-\rho_K) \phi_K m_K + \\
	+ \tau^n\sum_{\sigma\in\Sigma} F_{K,\sigma} \Big(\frac{\phi_L-\phi_K}{d_{\sigma}}\Big) m_{\sigma}d_{\sigma} +
	\sum_{K\in\mathcal{T}} E(\rho_K) m_K +\sum_{K\in\mathcal{T}} \lambda_K \rho_K m_K.
	\end{multline}
	The solution must satisfy the system of optimality conditions, namely the KKT conditions. Plugging the optimality condition w.r.t. $F_{K,\sigma}$, i.e. $F_{K,\sigma} = - (R_{\Sigma}(\boldsymbol{\rho}))_{\sigma} (\frac{\phi_L-\phi_K}{d_{\sigma}})$, in \eqref{eq:infsup_h} and considering that
	\[
	\sum_{\sigma\in\Sigma} (R_{\Sigma}(\rho^n))_{\sigma} \Big(\frac{\phi_L-\phi_K}{d_{\sigma}}\Big)^2 m_{\sigma}d_{\sigma}=\sum_{K\in\mathcal{T}} \rho_K \Big(R_{\mathcal{T}}\Big(\Big(\frac{\phi^n_L-\phi^n_K}{d_{\sigma}}\Big)^2 \Big)\Big)_Km_K,
	\]
	the optimality conditions reduce to the system
	\begin{equation}\label{eq:OC_h}
	\begin{cases}
	(\rho^n_K-\rho^{n-1}_K)m_K - \tau^n \sum_{\sigma\in\Sigma_K} (R_{\Sigma}(\rho^n))_{\sigma} (\frac{\phi^n_L-\phi^n_K}{d_{\sigma}}) m_{\sigma} = 0, \\
	(\phi^n_K-E'(\rho^n_K)-\lambda^n_K)m_K + \frac{\tau^n}{2} (R_{\mathcal{T}}((\frac{\phi^n_L-\phi^n_K}{d_{\sigma}})^2))_K m_K = 0, \\
	\rho^n_K \ge0, \; \lambda^n_K\le 0, \; \rho^n_K \lambda^n_K = 0,
	\end{cases}
	\quad \forall K \in \mathcal{T}.
	\end{equation}
	At each step $n$ of the discrete LJKO, the discrete density $(\rho^n_K)_{K \in \mathcal{T}}$ is completely defined by \eqref{eq:OC_h}.
	
	System \eqref{eq:OC_h} is not easy to solve, the major problem being the non-uniqueness of the multipliers $\boldsymbol{\lambda}$ and $\boldsymbol{\phi}$ whenever the density vanishes. 
	When upwinding is used for the reconstructed density, i.e. $\rho_{\sigma}=\rho_K$ if $\phi_L>\phi_K$, $\rho_{\sigma}=\rho_L$ otherwise, the Lagrange multiplier $\boldsymbol{\lambda}$ can be taken equal zero and disregarded \cite{LJKO_arXiv}. In our framework this is not possible and to avoid dealing explicitly with the positivity constraint we use an interior point method. The constraint is incoporated in the problem by adding to the functional a barrier function of the density which is convex and singular in zero. We use the logarithmic barrier $-\log(\rho)$. In this way the minimizer is automatically repulsed away from zero and the problem can be solved using the Newton scheme. The perturbation introduced by the barrier function can be tuned by multiplying it by a positive coefficient $\mu$.
	The perturbed version of problem \eqref{eq:infsup_h} for the $n-$th step of the discrete LJKO is
	\begin{multline}\label{eq:infsup_IPM}
	\inf_{(\boldsymbol{\rho},\boldsymbol{F})} \sup_{\boldsymbol{\phi}} \tau^n \sum_{\sigma\in\Sigma} \frac{(F_{\sigma})^2}{2(R_{\Sigma}(\boldsymbol{\rho}))_{\sigma}} m_{\sigma}d_{\sigma} + \sum_{K\in\mathcal{T}} (\rho^{n-1}_K-\rho_K) \phi_K m_K + \\
	+ \tau^n\sum_{\sigma\in\Sigma} F_{K,\sigma} \Big( \frac{\phi_L-\phi_K}{d_{\sigma}}\Big) m_{\sigma}d_{\sigma} +
	\sum_{K\in\mathcal{T}} E(\rho_K) m_K - \mu \sum_{K} \log(\rho_K) m_K,
	\end{multline}
	whose optimality conditions now are
	\begin{equation}\label{eq:OC_IPM}
	\begin{cases}
	(\rho^n_K-\rho^{n-1}_K)m_K - \tau^n \sum_{\sigma\in\Sigma_K} (R_{\Sigma}(\rho^n))_{\sigma} (\frac{\phi^n_L-\phi^n_K}{d_{\sigma}}) m_{\sigma} = 0, \\
	(\phi^n_K-E'(\rho^n_K)+s_K)m_K + \frac{\tau^n}{2} (R_{\mathcal{T}}((\frac{\phi^n_L-\phi^n_K}{d_{\sigma}})^2))_Km_K = 0, \\
	s_K \rho_K = \mu,
	\end{cases}
	\quad \forall K \in \mathcal{T},
	\end{equation}
	where the condition $F_{K,\sigma} = - (R_{\Sigma}(\boldsymbol{\rho}))_{\sigma} (\frac{\phi_L-\phi_K}{d_{\sigma}})$ has been substituted again.
	System \eqref{eq:OC_IPM} can be seen as a pertubation of \eqref{eq:OC_h}, where $\rho_K$ and $s_K=-\lambda_K$ are automatically forced to be positive and the orthogonality is relaxed. For small value of $\mu$ it provides an approximation of the solution $(\boldsymbol{\rho},\boldsymbol{\phi})$ to problem \eqref{eq:OC_h}. However, the smaller the parameter the more difficult it is to solve problem \eqref{eq:OC_IPM} with a Newton scheme.
	The idea is then to construct a sequence of solutions to problem \eqref{eq:OC_IPM} for a sequence of coefficients $\mu$ decreasing to zero, using the solution corresponding to the previous value of $\mu$ as starting point for the Newton scheme. In this way the solver approaches the solution to \eqref{eq:OC_h} from the interior of the region of feasibility: the density is always positive.
	\begin{algorithm}
		\SetAlgoLined
		Given the starting point $\boldsymbol{x}_0$ and the parameters $\mu_0>0,\theta\in(0,1),\varepsilon_0>0,\varepsilon_{\mu}>0$ \; 
		\While{$\delta_0>\varepsilon_0$}{
			$\mu = \theta \mu$ \;
			\While{$\delta_{\mu}>\varepsilon_{\mu}$}{
			compute Newton direction $\boldsymbol{d}$ for \eqref{eq:OC_IPM} and a step length $\alpha$\;
			update: $\boldsymbol{x} = \boldsymbol{x} + \alpha \boldsymbol{d}$ \;
			}
		}
		\caption{Interior point method}
		\label{alg:IPM}
	\end{algorithm}
	With reference to Algorithm \ref{alg:IPM}, $\varepsilon_0$ and $\varepsilon_{\mu}$ are the tolerances for the solution to \eqref{eq:OC_h} and \eqref{eq:OC_IPM} respectively, $\delta_0$ and $\delta_{\mu}$ denoting a norm of the residues of the two systems of optimality conditions. In practice, it is not necessary to find for each value of $\mu$ a precise solution, being interested only in the solution for $\mu=0$, and relatively big values can be used. Even doing only one Newton step, that is taking $\varepsilon_{\mu} = \infty$, can be sufficient and extremely effective.
	Moreover, the behavior of the solver strongly depends also on the initial value $\mu_0$ and the decay ratio $\theta \in (0,1)$, the difficulty to tune these parameters being its major drawback. We refer to \cite{Terlaky} and references therein for more details on interior point methods.
	
	As a final remark, note that solving the gradient flow with respect to an energy involving the entropy, i.e. $E(\rho) = \rho \log(\rho)$, automatically prevents the density from becoming negative. However, one cannot control the magnitude of the energy and therefore the interior point method, even if not strictly necessary, helps to get a more robust solver with respect to the Newton scheme. In fact, possible negative values for the density during the iterations of the algorithm could make it diverge, since the problem loses its convexity. The situation is similar when using the upwind technique to enforce the positivity. 

	\section{Numerical results}
	
	One of the most classical example of problems that exhibit a gradient flow structure is the Fokker-Planck equation:
	\begin{equation}
	\begin{cases}\label{eq:FokkerPlanck}
	\partial_t \rho = \Delta \rho + \nabla \cdot (\rho \nabla V) \quad &\text{in} \; \Omega\times [0,T], \\
	(\nabla \rho + \rho \nabla V) \cdot \boldsymbol{n}  = 0 \quad &\text{on} \; \partial\Omega\times [0,T],
	\end{cases}
	\end{equation}
	complemented with a positive initial condition, with  $V\in W^{1,\infty}(\Omega)$ a Lipschitz continuous exterior potential. 
	Equation \eqref{eq:FokkerPlanck} has been one of the first equations to be recasted as a gradient flow in the Wasserstein space with respect to the energy $\mathcal{E}(\rho) = \int_{\Omega} (\rho \log(\rho)+\rho V)\text{d} \boldsymbol{x}$ \cite{JKO98}. This example gives us the possibility to test the convergence of scheme \eqref{eq:ljko_h}. Consider indeed the density $\rho_s(\boldsymbol{x},t) = \exp(-(\pi^2+\frac{g^2}{4}) t+\frac{g}{2}x)(\pi \cos(\pi x)+\frac{g}{2}sin(\pi x))+\pi \exp(g(x-\frac{1}{2}))$, which is a solution to \eqref{eq:FokkerPlanck} in the domain $[0,1]^2\times[0,0.25]$ with potential $V(\boldsymbol{x})=-gx$. 
	Consider a sequence of meshes $\left(\mathcal{T}_m, \overline{\Sigma}_m, {(\boldsymbol{x}_K)}_{K\in\mathcal{T}_m}\right)$ with decreasing mesh size $h_m = h_{\mathcal{T}_m}$, and a sequence of decreasing time steps $\tau_m$ such that $(\frac{\tau_{m+1}}{\tau_m})=(\frac{h_{m+1}}{h_{m}})^2$. We solve problem \eqref{eq:FokkerPlanck} with scheme \eqref{eq:ljko_h} using this sequence of meshes and using as discrete initial condition $\rho_K^0 = \rho_s(\boldsymbol{x}_K,0)$. For each solution we compute the mesh-dependent $L^1((0,T);L^1(\Omega))$ error $\epsilon_m = \sum_{n} \tau_m \sum_{K \in \mathcal{T}_m} |\rho_K^n-\rho_s(\boldsymbol{x}_K,n \tau_m)| m_K$.
	In Table \ref{tab:errors} are listed the errors for each $m$ together with the convergence rate $\sqrt{\frac{\epsilon_{m-1}}{\epsilon_m}}$ for the three different weighted arithmetic averages. The scheme is first order accurate in time and second order accurate in space.
	\begin{table}
        \caption{Time-space convergence for the scheme.}
        \label{tab:errors}
        \begin{tabular}{cccccccc}
        \toprule
        $h_m$ & $\tau_m$ & $\epsilon_m^a$ & rate & $\epsilon_m^b$ & rate & $\epsilon_m^c$ & rate \\
        \midrule
		0.2986  &  0.0500 &  3.9382e-02  &    /      &  3.9526e-02  &    /  &  3.9157e-02  &     /  \\
		0.1493  &  0.0125 &  1.0345e-02  &  1.9286   &  1.0446e-02  &  1.9199 &  1.0246e-02  &  1.9342 \\
		0.0747  &  0.0031 &  2.6019e-03  &  1.9913   &  2.6367e-03  &  1.9861  &  2.5684e-03  &  1.9962 \\
		0.0373  &  0.0008 &  6.5090e-04  &  1.9990   & 6.6049e-04  &  1.9971 &  6.4170e-04  &  2.0009 \\
		0.0187  &  0.0002 &  1.6269e-04  &  2.0003   &  1.6519e-04  &  1.9994 &  1.6033e-04  &  2.0009 \\
        \bottomrule
        \end{tabular}
        \\ 
        $^a$ Weights $(\frac{1}{2},\frac{1}{2})$. $^b$ Weights $(\frac{d_L}{d_{\sigma}},\frac{d_K}{d_{\sigma}})$. $^c$ Weights $(\frac{d_K}{d_{\sigma}},\frac{d_L}{d_{\sigma}})$.
    \end{table}

\begin{figure}
		\centering
		\includegraphics[width=0.3\textwidth]{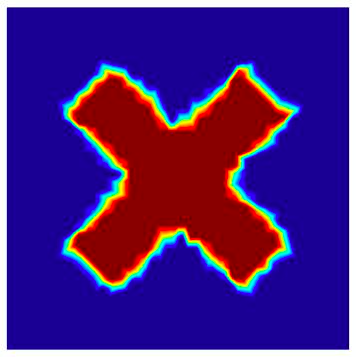}
		\includegraphics[width=0.3\textwidth]{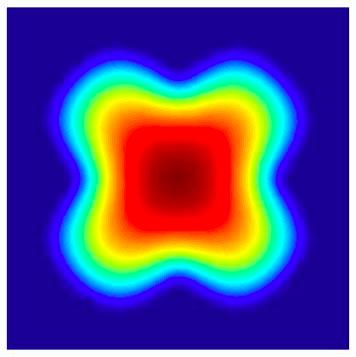}
		\includegraphics[width=0.3\textwidth]{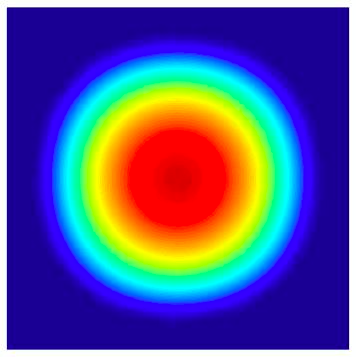}
		\includegraphics[height=0.3\textwidth]{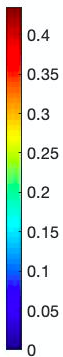}
		\caption{Convergence towards the Barenblatt solution ($\gamma=2$). Time steps $t=0$, $t=0.1$ and $t=0.7$.}
		\label{fig:porous_medium}
	\end{figure}

	As second application, we consider a gradient flow of an energy which is not singular in zero. On the domain $\Omega=[-1.5,1.5]^2$, for a time interval $[0,T]$, consider the porous medium equation,
\begin{equation*}
	\partial_t \rho = \Delta \rho^{\gamma} + \nabla \cdot (\rho \nabla V)
\end{equation*}
which has been proven in \cite{Otto01} to be a gradient flow in the Wasserstein space with respect to the energy $\mathcal{E}(\rho) = \int_{\Omega} \frac{1}{\gamma-1} \rho^{\gamma} + \rho V$, for a given $\gamma$ strictly greater than one. We consider the confining potential $V(\boldsymbol{x}) = \frac{1}{2} || \boldsymbol{x} ||^2_2$ which forces the density to concentrate at the origin.
In \eqref{fig:porous_medium} the evolution of an initial cross shaped density is shown for the case $\gamma=2$. As expected, the solution converges towards the Barenblatt profile $\rho^{\infty}(\boldsymbol{x}) = \max((\frac{M}{2\pi})^{\frac{\gamma-1}{\gamma}}- \frac{\gamma-1}{2\gamma}||\boldsymbol{x}||^2,0)^{\frac{1}{\gamma-1}}$, with $M$ being the total mass of the initial condition (Figure \ref{fig:porous_medium}).
	
	\begin{figure}
		\centering
		\includegraphics[height=0.3\textwidth]{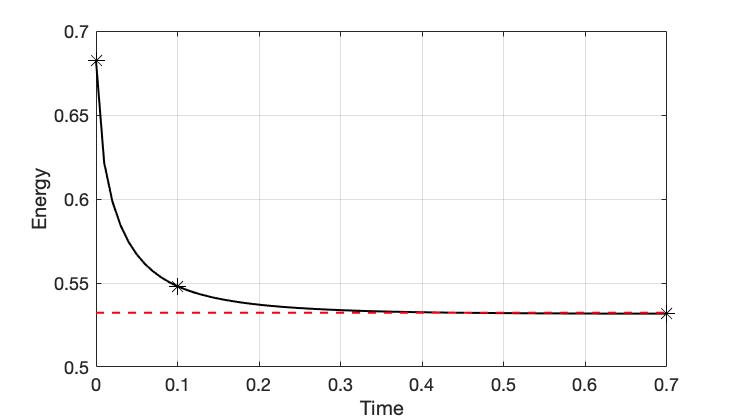} 
		\caption{Exponential decay profile of the discrete energy $\sum_{K\in\mathcal{T}} E(\rho_K) m_K$ (black), with the three values corresponding to Figure \ref{fig:porous_medium}, compared to the value of the energy for the Barenblatt equilibrium solution (red).}
		\label{fig:energy}
	\end{figure}

	
	\section*{Acknowledgements}
	The work of A. Natale was supported by the European Research Council (ERC project NORIA).
	G. Todeschi acknowledges that this project has received funding
	from the European Union’s Horizon 2020 research and innovation
	programme under the Marie Skłodowska-Curie grant agreement No 754362.
	\begin{center}
	\includegraphics[width=0.15\textwidth]{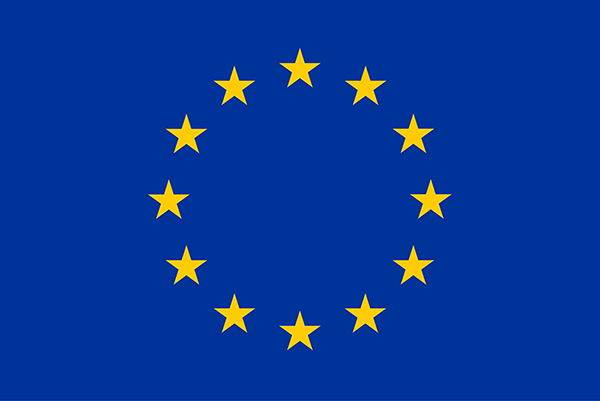}
	\end{center}

	\bibliography{biblio}

\end{document}